\documentclass{amsart}

 \theoremstyle{definition}
 
 \theoremstyle{remark}
 
 \numberwithin{equation}{subsection}
\begin{document}

\title[2-LOCAL DERIVATIONS ON VON NEUMANN ALGEBRAS ]
 {2-LOCAL DERIVATIONS ON SEMI-FINITE VON NEUMANN ALGEBRAS}

\author{ Shavkat Ayupov}

\address{National University of Uzbekistan, Institute of Mathematics,
Tashkent, Uzbekistan and the Abdus Salam International Centre for
Theoretical Physics (ICTP) Trieste, Italy}

\email{sh$_-$ayupov@mail.ru}

\author{ Farkhad Arzikulov}

\address{National University of Uzbekistan, Institute of Mathematics ,
Tashkent,  and Andizhan State University, Andizhan, Uzbekistan}

\email{arzikulovfn@rambler.ru}

%\thanks{This work was completed with the support of an Izaak
% Walton Killam Memorial Scholarship.}

%\thanks{The author was also supported in part by the Research
% Council of Slovenia.}

%\subjclass{Primary 46L57; Secondary 46L40}


\date{February 21, 2012.}

\dedicatory{}

%\commby{Daniel J. Rudolph}

%%% ----------------------------------------------------------------------

\begin{abstract}
In the present paper we prove that every 2-local
derivation on a semi-finite von Neumann algebra is a
derivation.
\end{abstract}

%%% ----------------------------------------------------------------------
\maketitle
%%% ----------------------------------------------------------------------

\section*{Introduction}

The present paper is devoted to 2-local derivations on von Neumann
algebras. Recall that a 2-local derivation is defined as follows:
given an algebra $A$, a map $\bigtriangleup : A \to A$ (not linear
in general) is called a 2-local derivation if for every $x$, $y\in
A$, there exists a derivation $D_{x,y} : A\to A$ such that
$\bigtriangleup(x)=D_{x,y}(x)$ and $\bigtriangleup(y)=D_{x,y}(y)$.

In 1997, P. \v{S}emrl \cite{S} introduced the notion of 2-local
derivations and described 2-local derivations on the algebra
$B(H)$ of all bounded linear operators on the infinite-dimensional
separable Hilbert space H. A similar description for the
finite-dimensional case appeared later in \cite{KK}. In the paper
\cite{LW} 2-local derivations have been described on matrix
algebras over finite-dimensional division rings.

In \cite{AK} the authors suggested a new technique and have
generalized the above mentioned results of \cite{S} and \cite{KK}
for arbitrary Hilbert spaces. Namely they considered 2-local
derivations on the algebra $B(H)$ of all linear bounded operators
on an arbitrary (no separability is assumed) Hilbert space $H$ and
proved that every 2-local derivation on $B(H)$ is a derivation.

In \cite{AA} we also suggested another technique and generalized
the above mentioned results of \cite{S}, \cite{KK} and \cite{AK}
for arbitrary von Neumann algebras of type I and proved
that every 2-local derivation on these algebras
is a derivation. In \cite{AKNA} (Theorem 3.4) a similar result was
proved for finite von Neumann algebras.

In the present paper we extended the above results and give a short proof of the theorem
for arbitrary semi-finite von Neumann algebras.

\section{Preliminaries}

Let $M$ be a von Neumann algebra.

{\it Definition.} A linear map $D : M\to M$ is called a
derivation, if $D(xy)=D(x)y+xD(y)$ for any two elements $x$, $y\in
M$.

A map $\Delta : M\to M$ is called a 2-local derivation, if for any
two elements $x$, $y\in M$ there exists a derivation $D_{x,y}:M\to
M$ such that $\Delta (x)=D_{x,y}(x)$, $\Delta (y)=D_{x,y}(y)$.

It is known that any derivation $D$ on a von Neumann algebra $M$
is an inner derivation, that is there exists an element $a\in M$
such that
$$
D(x)=ax-xa, x\in M.
$$
Therefore for a von Neumann algebra $M$ the above definition is
equivalent to the following one: A map $\Delta : M\to M$ is called
a 2-local derivation, if for any two elements $x$, $y\in M$ there
exists an element $a\in M$ such that $\Delta (x)=ax-xa$, $\Delta
(y)=ay-ya$.

Let $\mathcal{M}$ be a von Neumann algebra, $\bigtriangleup
:\mathcal{M}\to \mathcal{M}$ be a 2-local derivation. It easy to see
 that $\bigtriangleup$ is homogeneous. Indeed, for each $x\in
\mathcal{M}$, and for $\lambda \in {\mathbb C}$ there exists a
derivation $D_{x,\lambda x}$ such that
$\bigtriangleup(x)=D_{x,\lambda x}(x)$ and $\bigtriangleup(\lambda
x)=D_{x,\lambda x}(\lambda x)$. Then
$$
\bigtriangleup(\lambda x)=D_{x,\lambda x}(\lambda x)=\lambda
D_{x,\lambda x}(x) =\lambda \bigtriangleup(x).
$$
Hence, $\bigtriangleup$ is homogenous. Further, for each
$x\in \mathcal{M}$, there exists a derivation   $D_{x,x^2}$ such
that $\bigtriangleup(x)=D_{x,x^2}(x)$ and
$\bigtriangleup(x^2)=D_{x,x^2}(x^2)$. Then
$$
\bigtriangleup(x^2)=D_{x,x^2}(x^2)=D_{x,x^2}(x)x+xD_{x,x^2}(x)
=\bigtriangleup(x)x+x\bigtriangleup(x).
$$
A linear map satisfying the above identity is called a Jordan derivation.
In \cite{Bre} it is proved that any Jordan derivation on a
semi-prime algebra is a derivation. Since every von Neumann  algebra $\mathcal{M}$ is
semi-prime (i.e. $a\mathcal{M}a = \{0\}$ implies that $a = \{0\}$), in order to prove
that a 2-local derivation $\bigtriangleup :\mathcal{M}\to
\mathcal{M}$ is a derivation it is sufficient to show that the map
$\bigtriangleup :\mathcal{M}\to \mathcal{M}$ is additive.

\section{2-local derivations on semi-finite von Neumann algebras}

Let $\mathcal{M}$ be a semi-finite von Neumann algebra and
let $\tau$ be a faithful normal semi-finite trace on $\mathcal{M}$. Denote by $m_{\tau}$
the definition ideal of $\tau$, i.e the  set of all
elements $a\in \mathcal{M}$ such that $\tau(\vert a\vert)<\infty$. Then
$m_{\tau}$ is a $*$-algebra and, moreover $m_{\tau}$ is a two sided
ideal of $\mathcal{M}$  (see \cite{MT}, Definition 2.17).

It is clear that any derivation $D$ on $M$ maps the ideal $m_{\tau}$ into itself.
 Indeed, since $D$ is inner, i.e. $D(x)=ax-xa, x\in M$ for an appropriate $a\in M$,  we have that
  $D(x) = ax-xa \in m_{\tau}$ for all $x\in m_{\tau}$. Therefore any 2-local derivation on $M$ also maps $m_{\tau}$ into
itself.

{\bf Theorem.} {\it Let $\mathcal{M}$ be a semi-finite von Neumann algebra, and let
$\bigtriangleup :\mathcal{M}\to \mathcal{M}$
be a 2-local derivation. Then $\bigtriangleup$ is a derivation.}

{\it Proof.} Let $\triangle :\mathcal{M}\to \mathcal{M}$ be a 2-local derivation
and let $\tau$ be a faithful normal semi-finite trace on $\mathcal{M}$.
For each $x\in \mathcal{M}$ and $y\in m_{\tau}$ there exists a derivation
$D_{x,y}$ on $\mathcal{M}$ such that $\triangle(x)=D_{x,y}(x)$, $\triangle(y)=D_{x,y}(y)$.
Since every derivation on $\mathcal{M}$ is inner, there exists an element $a\in \mathcal{M}$
such that
$$
[a,xy]=D_{x,y}(xy)=D_{x,y}(x)y+xD_{x,y}(y)=\triangle(x)y+x\triangle(y)
$$
i.e.
$$
[a,xy]=\triangle(x)y+x\triangle(y).
$$
We have
$$
\vert \tau(axy)\vert<\infty.
$$
Since $m_{\tau}$ is an ideal and $y\in m_{\tau}$, the elements $axy$, $xy$, $xya$ and $\triangle(y)$ also belong to $m_{\tau}$ and hence we have
$$
\tau(axy)=\tau(a(xy))=\tau((xy)a)=\tau(xya).
$$
Thus
$$
0=\tau(axy-xya)=\tau([a,xy])=\tau(\triangle(x)y+x\triangle(y)),
$$
i.e.
$$
\tau(\triangle(x)y)=-\tau(x\triangle(y)).
$$
For arbitrary $u$, $v\in \mathcal{M}$ and $w\in m_{\tau}$ set $x=u+v$, $y=w$. Then $\triangle(w)\in m_{\tau}$ and
$$
\tau(\triangle(u+v)w)=-\tau((u+v)\triangle(w))=
$$
$$
-\tau(u\triangle(w))-\tau(v\triangle(w))=\tau(\triangle(u)w)+\tau(\triangle(v)w)=
$$
$$
\tau((\triangle(u)+\triangle(v))w),
$$
and so
$$
\tau((\triangle(u+v)-\triangle(u)+\triangle(v))w)=0,
$$
for all $u$, $v\in \mathcal{M}$ and $w\in m_{\tau}$. Denote $b=\triangle(u+v)-\triangle(u)+\triangle(v)$.
Then
$$
\tau(bw)=0 \,\,\,\, \forall w\in m_{\tau}          \,\,\,\, (1).
$$
Now take a monotone increasing net $\{e_{\alpha}\}_{\alpha}$ of projections in
$m_{\tau}$ such that $e_{\alpha}\uparrow 1$ in $\mathcal{M}$.
Then $\{e_{\alpha}b^*\}_{\alpha} \subset m_{\tau}$. Hence  (1) implies
$$
\tau(be_{\alpha}b^*)=0 \,\,\,\, \forall \alpha.
$$
At the same time $be_{\alpha}b^*\uparrow bb^*$ in $\mathcal{M}$.
Since the trace $\tau$ is normal we have
$$
\tau(be_{\alpha}b^*)\uparrow \tau(bb^*),
$$
i.e.  $\tau(bb^*)=0$.  The trace $\tau$ is faithful so this implies that  $bb^*=0$, i.e. $b=0$.
Therefore
$$
\triangle(u+v)=\triangle(u)+\triangle(v), u,v\in \mathcal{M},
$$
i.e. $\triangle$ is an additive map on $\mathcal{M}$. As it was mentioned in "Preliminaries" this implies
that  $\triangle$
is a derivation on $\mathcal{M}$. The proof is complete.
$\triangleright$

\end{document}